\documentclass{article}

\newtheorem{fed}{\textbf{Definition}}[section]
\newtheorem{thm}[fed]{\textbf{Theorem}}
\newtheorem{lemma}[fed]{\textbf{Lemma}}

\newtheorem{rem}[fed]{\textbf{Remark}}
\newtheorem{prop}[fed]{\textbf{Proposition}}

\usepackage{amssymb,bbm,graphicx,epsfig,psfrag,epic,eepic,latexsym}
\usepackage{amsmath}
\usepackage{mathrsfs} 
\usepackage[dvips]{color}       

\begin{document}
\title{Vortices on the cylinder}
\author{Urs Frauenfelder\footnote{
Supported by JSPS}}
\maketitle

\begin{abstract}
We apply the finite dimensional approximation techniques of
Furuta, Kronheimer, and Manolescu to give a new proof
of a result of Jaffe and Taubes.
\end{abstract}
\tableofcontents

\newpage

\section[Introduction]{Introduction}

The vortex equations are the absolute minima of the Yang-Mills-Higgs 
functional. For a unitary line bundle $L$ over a Riemann surface 
let $\mathfrak{A}$ be the space of unitary connections of $L$ and
$\Omega(L)$ be the space of smooth sections of $L$ then the Yang-Mills-Higgs 
functional $YMH \colon \mathfrak{A}\times \Omega(L)\to \mathbb{R}$
is defined as
$$YMH(A,v)=\int_{\Sigma}\big(|F_A|^2+|\nabla_A v|^2
+\frac{1}{4}(1-|v|^2)^2\big) dvol.$$
Its absolute minima
\begin{eqnarray}\label{vortexeq}
\bar{\partial}_A v=0\\ \nonumber
*F_A=\frac{1}{2}\big(|v|^2-1\big)
\end{eqnarray}
are the vortex equations. The Yang-Mills-Higgs functional is invariant
under the action of the gauge group $\mathcal{G}=C^\infty(\Sigma,S^1)$
and hence so are the vortex equations. 

For the case $\Sigma=\mathbb{C}$ the moduli spaces for the vortex equations
were completely described by Jaffe and Taubes, see \cite{jaffe-taubes}.
For solutions $(A,v)$ of the vortex equations on $\mathbb{C}$ satisfying
an appropriate decay condition at infinity, it turns out that
the vortex number
$$N=\frac{1}{2\pi}\int_{\mathbb{C}}F_A$$
is an integer. Jaffe and Taubes proved that the moduli space of vortices
with vortex number $N$ modulo gauge is given by the $N$-fold symmetric product 
$$\{(\ref{vortexeq}): \mathrm{vortex\,\,number}=N\}/\mathcal{G}
\cong
S^N\mathbb{C} \cong \mathbb{C}^N.$$
The case for compact Riemann surfaces $\Sigma$ was studied by Bradlow and
Garcia-Prada, see \cite{bradlow, garcia-prada}. In the compact case the
vortex number is 
$$N=\langle c_1(L),[\Sigma]\rangle$$
and the moduli space was determined by Bradlow and Garcia-Prada
$$\{(\ref{vortexeq}): \mathrm{vortex\,\,number}=N\}/\mathcal{G}
\cong
\left\{\begin{array}{cc}
S^N \Sigma & N<vol(\Sigma)/4\pi\\
\emptyset & N > vol(\Sigma)/4\pi.
\end{array}\right.
$$
If $N=vol(\Sigma)/4\pi$ then solutions of (\ref{vortexeq}) necessarily
satisfy $v\equiv 0$.
In this article we consider the case where $\Sigma=\mathcal{Z}$ is
the cylinder. We prove \\ \\
\textbf{Theorem A}  \emph{The moduli space $N$-vortices 
on the cylinder modulo gauge is 
$S^N \mathcal{Z}$.}
\\ \\
We do not claim originality for this theorem since the methods of
Jaffe and Taubes for the complex plane could also be used
to determine the vortices on the cylinder. However, we will present 
in this paper a new approach for proving existence of PDE's
by using the finite dimensional approximation 
techniques of Furuta, Kronheimer, and
Manolescu \cite{furuta, kronheimer-manolescu, manolescu}. 
\\
The idea of this new method is the following. Solutions
of the vortex equations on the cylinder can be
interpreted as 
flow lines of an action functional $\mathcal{A}$
defined on an infinite dimensional space $\mathscr{L}$.
We consider a finite dimensional approximation
$L \subset \mathscr{L}$ and homotop the flow lines of
$\mathcal{A}$ to the flow lines of $\mathcal{A}|_L$.
Since $L$ is finite dimensional the flow lines of the
restricted action functional are solutions of an ODE.
This enables us to translate the question of existence of
a PDE to the question of existence of an ODE.\\
However, to prove existence of finite energy Morse flow lines on a
noncompact manifold is still a hard task. To do that
we will take advantage of the fact that the restriction of
our action functional to the finite dimensional approximation
has the form of a Lagrange multiplier functional. It is well
known from basic calculus that critical points of 
a function under a constraint can be found by considering
the Lagrange multiplier functional. However, the Morse flow lines
of the Lagrange multiplier functional and the Morse flow lines
of the function restricted to the constraint are in
general quite different. We will develop a theory which shows
how they can be homotoped to each other. This theory allows us
to translate the question of existence of Morse flow lines
on a noncompact manifold to the question of existence 
of Morse flow lines on a compact manifold.  
\\ \\
This paper is organized as follows. In Section~\ref{A}
we prove Theorem A using finite dimensional approximation modulo
the theory of Morse functions with Lagrange multipliers. 
In Section~\ref{Further} we discuss further examples were our
methods could be applied. In the Appendix we discuss the Theory
of Morse functions with Lagrange multipliers.
\\ \\
\textbf{Acknowledgements: }I would like to thank Kaoru Ono
for pointing out to me that the finite dimensional approximation 
techniques of Furuta, Kronheimer, and Manolescu can be applied
to the vortex equations on the cylinder. The final part of this
paper was written at FIM of ETH Z\"urich. I wish to thank FIM for
its kind hospitality.

\newpage

\section[Proof of Theorem A]{Proof of Theorem A}\label{A} 

\subsection[The gradient equation]{The gradient equation}

The standard circle action on $\mathbb{C}$ given by
$$z \mapsto e^{i\theta}z, \quad e^{i\theta} \in S^1$$
is Hamiltonian with respect to the standard symplectic structure
$\omega=dx \wedge dy$ on $\mathbb{C}$. A moment map for the action
is given by
$$\mu(z)=-\frac{i}{2}|z|^2+\frac{i}{2}
\in i\mathbb{R}=\mathrm{Lie}(S^1).$$ 
We consider the loop space
$$\mathscr{L}=C^\infty(S^1,\mathbb{C}\times i\mathbb{R})$$
and define the action functional $\mathcal{A}\colon
\mathscr{L} \to \mathbb{R}$ by
\begin{equation}\label{action}
\mathcal{A}(v,\eta)=\int_0^1\lambda(v)(\partial_t v)+
\int_0^1\langle \mu(v),\eta\rangle dt
\end{equation}
where $\lambda=y dx$ is the Liouville one-form on
$\mathbb{C}$ satisfying $d \lambda=-\omega$. The first term 
in (\ref{action}) is just Floer's action functional
on the loop space $C^\infty(S^1,\mathbb{C})$, 
$$\mathcal{A}_{fl}(v)=\int_0^1 \lambda(v)(\partial_t v)$$
where the second term may be thought of as a Lagrange multiplier to
the constraint $\mu^{-1}(0)$. 

The gauge group
$$\mathcal{H}=C^\infty(S^1,S^1)$$
acts on $(v,\eta) \in \mathscr{L}$ by
$$h_*(v,\eta)=(hv,\eta-h^{-1}\partial_t h),\quad
h \in \mathcal{H}.$$
The differential of the action functional $d \mathcal{A}$ and
the $L^2$-metric $g_{L^2}$ on $\mathscr{L}$ are invariant under
the gauge action and hence also the
gradient flow lines of $\nabla_{g_{L^2}}\mathcal{A}$ which are
solutions $(v,\eta)\in C^\infty(\mathbb{R}\times S^1, \mathbb{C}
\times i\mathbb{R})$ of the following PDE
\begin{eqnarray}\label{gradeq}
\partial_s v+i \partial_t v+i \eta v=0\\ \nonumber
\partial_s \eta+\mu(v)=0.
\end{eqnarray}
Solutions of (\ref{gradeq}) are solutions of (\ref{vortexeq}) in
radial gauge. One has a natural bijection
$$\{(\ref{vortexeq})\}/\mathcal{G}\cong 
\{(\ref{gradeq})\}/\mathcal{H}$$
by setting $A=\eta dt+\zeta ds$ and using a gauge transformation
$g \in \mathcal{G}$ such that $g_*\zeta=0$.
\begin{rem}\label{sve}
There is a straightforward generalization of the
action functional (\ref{action}) to general Hamiltonian group actions
on symplectic manifolds. The gradient flow lines of these action 
functionals are the symplectic vortex equations 
\cite{cieliebak-gaio-mundet-salamon, cieliebak-gaio-salamon, frauenfelder}
which reduce in
the special case of a circle action on 1-dimensional complex space
to the classical vortex equations. 
\end{rem}
Since the Marsden-Weinstein quotient $\mu^{-1}(0)/S^1$ is just a point,
the critical manifold of $\mathcal{A}$ is homeomorphic to the gauge
group $\mathcal{H}$. In particular,
$$\pi_0(\mathrm{crit}(\mathcal{A}))\cong
\pi_0(\mathcal{H})\cong \pi_1(S^1) \cong \mathbb{Z}$$
which enables us to recover the vortex number in this setting. 
The vortex number of a solution of (\ref{gradeq}) 
is proportional to the energy via
$$N=\frac{1}{\pi}\lim_{s \to \infty}\Big(\mathcal{A}\big((v,\eta)(s,\cdot)\big)
-\mathcal{A}\big((v,\eta)(-s,\cdot)\big)\Big)=\frac{1}{\pi}E(v,\eta)$$
We denote by $\mathfrak{V}^N$ the moduli space of $N$-vortices on
the cylinder, i.e. gradient flow lines of 
$\nabla_{g_{L_2}}\mathcal{A}$ modulo gauge which converge at the ends
to connected components of $\mathrm{crit}(\mathcal{A})$ of
difference $N \in \mathbb{Z}$. The following theorems follow from the
results in \cite{cieliebak-gaio-mundet-salamon}.
\begin{thm}[Regularity]
The moduli spaces $\mathfrak{V}^N$ are smooth manifolds of dimension
$2N$.
\end{thm}

\begin{thm}[Compactness moduli breaking] Let
$(v^\nu, \eta^\nu)$ be a sequence of $N$-vortices. Then there exists
a subsequence $\nu_j$, a sequence of gauge transformations
$h_j \in \mathcal{H}$, $N_i$-vortices $(v_i,\eta_i)$ for
$1 \leq i \leq \ell$, and sequences of real numbers $S^j_i$ such that
the timeshifted vortices converge uniformly in the
$C^\infty_{loc}$-topology
$$(h_j)_*(v^{\nu_j},\eta^{\nu_j})(\cdot, \cdot+S^j_i)
\longrightarrow_{j \to \infty} (v_i,\eta_i)$$
different timeshifts diverge
$$\lim_{j \to \infty}|S^j_i-S^j_{i'}|=\infty, \quad i \neq i'$$
and the total vortex number is preserved
$$\sum_{i=1}^\ell N_i=N.$$
\end{thm}

\subsection[Finite dimensional approximation]{Finite dimensional approximation}

The gauge group decomposes
$$\mathcal{H}=\mathcal{H}_0 \oplus S^1 \oplus \mathbb{Z}.$$
where the infinite dimensional contractible group
$\mathcal{H}_0$ is given by
$$\mathcal{H}_0=\Bigg\{g=\exp(\xi) \in \mathcal{H}: 
\xi \in C^\infty(S^1, i\mathbb{R}),\,\,\int_0^1 \xi dt=0\Bigg\}.$$
Following Manolescu \cite{manolescu} we get rid of 
$\mathcal{H}_0$ by projecting the gradient equations to the Coulomb section
in $\mathscr{L}$. To see how this works, observe that 
$\mathcal{H}_0$ acts freely on $\mathscr{L}$ and 
since our gauge group is abelian we can put 
each $\eta \in C^\infty(S^1,i\mathbb{R})$ into global Coulomb gauge on
the circle, namely there exists a unique
$h_\eta \in \mathcal{H}_0$ such that 
$$0=d^*\big((h_\eta)_*\eta\big)=-\partial_t\big((h_\eta)_*\eta\big).$$
Hence we may think of 
$$\mathscr{L}_c=C^\infty(S^1,\mathbb{C})\times i\mathbb{R}$$
as a section in the principal $\mathcal{H}_0$-bundle 
$\mathscr{L}$, or more precisely, we have a commutative diagram
\\ \\
\setlength{\unitlength}{1.5cm}
\begin{picture}(5,2)\thicklines
 \put(2,0){$\mathscr{L}_c$}
 \put(2.3,0.3){\vector(2,1){2}}
 \put(3.2,0.9){$\iota$}
 \put(2.5,0.1){\vector(1,0){1.7}}
 \put(3.4,0,2){$c$}
 \put(4.3,0){$\mathscr{L}/\mathcal{H}_0$}
 \put(4.4,1.2){$\mathscr{L}$}
 \put(4.5,1.1){\vector(0,-1){0.8}}
 \put(4.6,0.6){$\pi$}
\end{picture}
\\ \\
where $\iota$ denotes the canonical inclusion and
$c$ denotes the bijection which is induced by Coulomb gauge.

The $L^2$-metric $g_{L^2}$ on $\mathscr{L}$ induces two natural metrics
on $\mathscr{L}_c$
$$g_0=\iota^*g_{L^2}, \quad g_1=c^*[g_{L^2}]$$
where $[g_{L^2}]$ denotes the quotient metric of the $L^2$-metric
on $\mathscr{L}/\mathcal{H}_0$.

Abbreviate 
$$\mathcal{A}_c=\mathcal{A}|_{\mathscr{L}_c}.$$
Then $\mathcal{H}_0$-gauge equivalence classes of 
flow lines of $\nabla_{g_{L^2}}\mathcal{A}$ are in natural one-to-one
correspondence with flow lines of $\nabla_{g_1}\mathcal{A}_c$ by
projection. The importance of $g_0$ lies in the fact that flow lines 
of $\nabla_{g_0}\mathcal{A}_c$ are
contained in natural finite dimensional subspaces of the infinite
dimensional space $\mathscr{L}_c$. For
integers $\mu \leq \nu$ consider the Fourierapproximations
$$L^\nu_{\mu}=\{z=\sum_{j=\mu}^\nu z_j e^{2\pi i j}: z_j \in \mathbb{C}\}$$
of the loop space $C^\infty(S^1,\mathbb{C})$. The metric $g_0$
is just the product of the $L^2$-metric on $C^\infty(S^1,\mathbb{C})$
and the metric induced from the inner product on $i\mathbb{R}$ and hence
$$\nabla_{g_0}\mathcal{A}_c(z,\eta) \in L_{\mu}^\nu \times i\mathbb{R}
\subset \mathscr{L}_c, \quad (z,\eta) \in L^\nu_\mu\times i\mathbb{R}.$$ 
It follows that flow lines of 
$\mathcal{A}_c|_{L^\nu_\mu \times i\mathbb{R}}$ are actually flow lines
of $\mathcal{A}_c$. Moreover, critical points of $\mathcal{A}_c$
are tuples
$$(v_0 e^{2\pi i m t}, 2\pi m i), \quad |v_0|=1,\,\, m \in \mathbb{Z}$$
and hence for every pair of critical points of
$\mathcal{A}_c$ there exists a finite dimensional approximation as above
which contains both of them. 

The following proposition shows that for
each finite energy flow line of $\mathcal{A}_c$ one can find 
a finite dimensional approximation such that the flow line is entirely
contained in it. For a set of 
flow lines which converge at both ends to the same
critical points there can be found a finite dimensional approximation which
contains the whole set simultanuously. However note, that since there
are infinitely many critical points of $\mathcal{A}_c$ there is no
finite dimensional approximation in which all finite energy flow lines lie 
simultanuously.    

\begin{prop}\label{finiteapp}
Assume that $(v,\eta) \in C^\infty(\mathbb{R}\times S^1,\mathbb{C})
\times C^\infty(\mathbb{R},i\mathbb{R})$ is a gradient flow line of
$\nabla_{g_0} \mathcal{A}_c$ such that
$$\lim_{s \to \pm \infty}(v,\eta)(s,t)=(v_\pm e^{2\pi i m_\pm t}, 
2 \pi m_\pm i), \quad |v_\pm|=1,\,\,m_\pm \in \mathbb{Z}$$
where the limit is uniformly with respect to the $C^\infty$-topology. 
Then $(v,\eta)(s,\cdot)$ is contained in 
$L^{-m_+}_{-m_-}\times i\mathbb{R}$. 
\end{prop}
\textbf{Proof: }Abbreviate $\bar{\mu}(v) \in C^\infty(\mathbb{R},i\mathbb{R})$
by
$$\bar{\mu}(v)(s)=\int_0^1 \mu(v(s,t))dt, \quad s \in \mathbb{R}.$$
A gradient flow line of $\nabla_{g_0}\mathcal{A}_c$
is a solution of the following PDE
\begin{eqnarray}\label{grad0}
\partial_s v+i\partial_t v+i\eta v=0\\ \nonumber
\partial_s \eta+\bar{\mu}(v)=0.
\end{eqnarray}
Plugging in the Fourierexpansion
$$v(s,t)=\sum_{m=-\infty}^\infty v_j(s)e^{2\pi i jt}$$
into the first equation of (\ref{grad0}) we obtain
for each Fouriercoefficient the ODE
\begin{equation}\label{fouriercof}
\partial_s v_m(s)+\big(i\eta(s)-2\pi m\big)v_m(s)=0. 
\end{equation} 
Using (\ref{fouriercof}), $\lim_{s\to \pm \infty}\partial_s v_m(s)=0$,
and the asymptotic behaviour of $\eta(s)$, we conclude that
$v_m$ vanishes identically unless $m$ is contained in 
$\{-m_-,\ldots,-m_+\}$. This proves the proposition.
\hfill $\square$
\\ \\
In order to homotop the PDE (\ref{gradeq}) to an ODE it remains 
in view of the proposition above
to find a homotopy between $g_0$ and $g_1$. This homotopy has to be
compact, i.e. the moduli spaces of finite energy flow lines should be
compact modulo breaking and modulo the remaining action of
the noncompact group $\mathbb{Z}$. Moreover, we require the homotopy to be
equivariant with respect to the following torus action. There is the
circle action on the target manifold $\mathbb{C}$ and there is a further
circle action on the domain $S^1$ given by rotating the circle. 
The two actions commute on $\mathscr{L}_c$ and lead to an action
of the two torus $T^2=S^1 \times S^1$
on $\mathscr{L}_c$. Note that the action functional $\mathcal{A}_c$, 
and the metrics $g_0$ and $g_1$ are $T^2$-invariant.  
\begin{thm}\label{comhop}
There exists a continuous family of $T^2$-invariant metrics 
$g_r$ for $r \in [0,1]$ on $\mathscr{L}_c$ with the following property.
Assume that for $\nu \in \mathbb{N}$ there exists a sequence
of flow lines $(v^\nu,\eta^\nu)$ of $\nabla_{g_{r^\nu}}\mathcal{A}_c$
for $r^\nu \in [0,1]$
whose energy is uniformly bounded, i.e. there exists
a constant $c>0$ such that for all $\nu$ it holds
$$E(v^\nu,\eta^\nu) \leq c.$$
Then there exists a subsequence $\nu_j$, a sequence of gauge 
transformations $h_j \in \mathbb{Z}$, flow lines
$(v_i,\eta_i)$ for $1 \leq i \leq \ell$ of
$\nabla_{g_{r_\infty}}\mathcal{A}_c$ for $r_\infty \in [0,1]$,
and sequences of real numbers $S^j_i$ such that the timeshifted
vortices converge uniformly in the $C^\infty_{loc}$-topology
$$(h_j)_*(v^{\nu_j},\eta^{\nu_j})(\cdot, \cdot+S^j_i)
\longrightarrow_{j \to \infty} (v_i,\eta_i)$$
different timeshifts diverge
$$\lim_{j \to \infty}|S^j_i-S^j_{i'}|=\infty, \quad i \neq i'$$
and the total energy is preserved
$$\lim_{j \to \infty}E(v^{\nu_j},\eta^{\nu_j})=\sum_{i=1}^\ell
E(v_i,\eta_i)$$
\end{thm}
\textbf{Proof:} 
We first construct a $T^2$-invariant homotopy between $g_0$ and
$g_1$. 
In order to do that, observe that the geometric reason 
that $g_0$ and $g_1$ are different lies in the fact that
the infinitesimal gauge action of $\mathcal{H}_0$ is
not orthogonal to the Coulomb section $\mathscr{L}_c$ 
with respect to the $L^2$-metric on $\mathscr{L}$.
To construct the homotopy we consider a family of $\mathcal{H}_0$ 
actions on $\mathscr{L}$ such that $\mathscr{L}_c$ is a section
for the whole family of actions but that for the final action
the Coulomb section gets orthogonal to the infinitesimal gauge
action.
 
Taking advantage of the contractibility of
the gauge group $\mathcal{H}_0$ we define for $r \in [0,1]$ 
and $h=\exp(\xi) \in \mathcal{H}_0$ the $h_{*_r}$ action on
$(v,\eta) \in \mathscr{L}$ by
$$h_{*_r}(v,\eta)=(\exp(r\xi)v,\eta-h^{-1}\partial_t h).$$
The deformed actions of $\mathcal{H}_0$ are still free on
$\mathscr{L}$ and $\mathscr{L}_c$ is a simultanuous section for
the whole family of actions. For each $r \in [0,1]$ we have
a commutative diagram
\\ \\
\setlength{\unitlength}{1.5cm}
\begin{picture}(5,2)\thicklines
 \put(2,0){$\mathscr{L}_c$}
 \put(2.3,0.3){\vector(2,1){2}}
 \put(3.2,0.9){$\iota$}
 \put(2.5,0.1){\vector(1,0){1.7}}
 \put(3.4,0,2){$c_r$}
 \put(4.3,0){$\mathscr{L}/_r\mathcal{H}_0$}
 \put(4.4,1.2){$\mathscr{L}$}
 \put(4.5,1.1){\vector(0,-1){0.8}}
 \put(4.6,0.6){$\pi_r$}
\end{picture}
\\ \\
where $\mathscr{L}/_r \mathcal{H}_0$ denotes the quotient of
$\mathscr{L}$ under the $r$-action of $\mathcal{H}_0$,
$\pi_r$ denotes the according canonical projection, and
$c_r$ refers to the Coulomb gauge of the $r$-action. The 
$L^2$-metric on $\mathscr{L}$ is simultanuously
$\mathcal{H}_0$ invariant for the whole family of actions and
hence induces for every $r \in [0,1]$ a quotient metric
$[g_{L^2}]_r$ on $\mathscr{L}/_r \mathcal{H}_0$. We define
$$g_r=c_r^*[g_{L^2}]_r.$$
It is easy to check that $g_r$ are $T^2$-invariant for every
$r \in [0,1]$. Moreover, for $r=0$ the Coulomb section is
orthogonal to the $0$-action of $\mathcal{H}_0$ and
hence $g_0$ defined in this way agrees with the previous definition
of $g_0$.

The gradient flow lines of $\mathcal{A}_c$ with respect to the
metric $g_r$ are solutions of the following problem
\begin{eqnarray}\label{homeq}
\partial_s v+\xi_v v+i\partial_t v+i \eta v=0\\ \nonumber
\partial_s \eta+\bar{\mu}(v)=0
\end{eqnarray}
where $\xi_v \in C^\infty(\mathbb{R}\times S^1,i\mathbb{R})$
which is determined for every $s \in \mathbb{R}$ by the conditions
$$\partial_t \xi_v(s,\cdot)=r^2\big(\mu(v(s,\cdot))-\bar{\mu}(v)(s)\big),\quad
\int_0^1 \xi_v(s,t)dt=0.$$
The main difficulty for proving the compactness statement in 
Theorem~\ref{comhop} is to obtain a uniform $L^\infty$-estimate
independent of $r \in [0,1]$ for all finite energy solutions 
of (\ref{homeq}). This provides the following lemma.

\begin{lemma} Let $(v,\eta)\in C^\infty(S^1 \times \mathbb{R},\mathbb{C})
\times C^\infty(\mathbb{R},i\mathbb{R})$ be a finite energy solution
of (\ref{homeq}) for
$r \in [0,1]$. Then there exists a constant $c<\infty$ independent
of $r$ such that $||v||_\infty<c$
\end{lemma}
\textbf{Proof: }Define $u(s)=
\frac{1}{2}\int_0^1 |v(s,t)|^2dt$ for $s \in \mathbb{R}$.
\\ \\
\textbf{Step 1: }\emph{$u(s) \leq 1/2$ for every $s \in \mathbb{R}$.}
\\ \\
Using a computation similar to the one in the proof of
\cite[Proposition 3.5]{cieliebak-gaio-salamon}
we estimate
\begin{eqnarray*}
\partial^2_s u
&=&\int_0^1\big(|\partial_s v+\xi_v v|^2
+|\partial_t v+\eta v|^2\big)dt+\\ 
& &2\int_0^1\big\langle \mu(v),(1-r^2)\bar{\mu}(v)+r^2
\mu(v)-i/2\big\rangle dt\\ 
&\geq&2\langle \bar{\mu}(v),\bar{\mu}(v)-i/2\rangle\\ 
&=&2\langle u,u+1/2\rangle\\ 
&\geq&2u(u-1/2).
\end{eqnarray*}
Hence if $u(s_0)>1/2$ for $s_0 \in \mathbb{R}$, then 
$u$ cannot have a local maximum at $s_0$. However the finite
energy assumption implies that $\lim_{s \to \pm \infty}u(s)=1/2$
which proves Step 1. 
\\ \\
\textbf{Step 2: }\emph{There exists a constant $c_1$
and a gauge transformation
$h \in \mathbb{Z}$ such that
$||\partial_t (h_*v)||_\infty \leq c_1||h_*v||_\infty^2
=c_1||v||^2_\infty.$}
\\ \\
Fix some integer $n>3$ and consider the finite cylinder
$\mathcal{Z}_n=S^1 \times [-n,n]$. 
It follows from Step 1 that
\begin{equation}\label{e1}
||v||_{L^2(\mathcal{Z}_n)}=\mathcal{O}(1).
\end{equation}
After a gauge transformation we may assume without
loss of generality that
$$||\eta(0)||=\mathcal{O}(1).$$
Using the second equation in (\ref{homeq}) 
and Step 1 we conclude that
\begin{equation}\label{e2}
||\eta||_{L^\infty(\mathcal{Z}_n)}=\mathcal{O}(1).
\end{equation}
The definition of $\xi_v$ together with Step 1 implies
that
\begin{equation}\label{e3}
||\xi_v||_{L^\infty(\mathcal{Z}_n)}=\mathcal{O}(1).
\end{equation}
Combining (\ref{e1}),(\ref{e2}), and (\ref{e3}) and using the
first equation in (\ref{homeq}) we conclude that
$$||\bar{\partial}v||_{L^2(\mathcal{Z}_n)}=\mathcal{O}(1)$$
from which we deduce using (\ref{e1}) and elliptic regularity
for the Cauchy-Riemann operator
\begin{equation}\label{e4}
||v||_{W^{1,2}(\mathcal{Z}_{n-1})}=\mathcal{O}(1).
\end{equation}
It follows from Sobolev's embedding theorem that for every $p<\infty$
we have
\begin{equation}\label{e41}
||v||_{L^p(\mathcal{Z}_{n-1})}=\mathcal{O}_p(1)
\end{equation}
from which we deduce analogously as before
\begin{equation}\label{e42}
 ||v||_{W^{1,p}(\mathcal{Z}_{n-2})}=\mathcal{O}_p(1).
\end{equation}
Using (\ref{e42}) and the definition of $\xi_v$ we
conclude
\begin{equation}\label{e5}
||\xi_v||_{W^{1,p}(\mathcal{Z}_{n-2})}=\mathcal{O}_p(||v||_\infty).
\end{equation}
The Laplacian of $v$ satisfies the equation
\begin{eqnarray}\label{lap}
\Delta v&=&i (\partial_t \xi_v) v+i \xi_v(\partial_tv)
-(\partial_s \xi_v) v-\xi_v (\partial_s v)\\ \nonumber
& &-i(\partial_s \eta) v-i\eta(\partial_s v)
-\eta (\partial_t v).
\end{eqnarray}
Using (\ref{e5}) and (\ref{e42}) we conclude from (\ref{lap}) that
$$||\Delta v||_{L^p(\mathcal{Z}_{n-2})}=\mathcal{O}_p(||v||_\infty^2)$$
from which we conclude by elliptic regularity for the Laplace operator
and (\ref{e41})
\begin{equation}\label{e6}
||v||_{W^{2,p}(\mathcal{Z}_{n-3})}=\mathcal{O}_p(||v||_\infty^2).
\end{equation}
Step 2 follows now from (\ref{e6}) and the Sobolev embedding theorem. 
\\ \\
\textbf{Step 3: }\emph{Proof of the lemma.}
\\ \\
Abbreviate $v_s=v(s,\cdot)$ and let $||v_s||_p$ be the $L^p$-norm on
the circle. It follows from Step 1 and Step 2 that there exist
constants $c_0$ and $c_1$ such that
\begin{equation}\label{ee1}
||v_s||_2 \leq c_0, \quad ||\partial_t v_s||_\infty \leq c_1 ||v_s||_\infty^2.
\end{equation}
We may assume without loss of generality that
$$|v_s(0)|=||v_s||_\infty.$$
We then estimate for $t \in S^1=\mathbb{R}/\mathbb{Z}$
using the second inequality in (\ref{ee1})
\begin{equation}\label{ee2}
|v_s(t)| \geq ||v_s||_\infty-c_1||v_s||_\infty^2|t|.
\end{equation}
Hence
\begin{eqnarray*}
||v_s||_2 &\geq& 
\Bigg(2\int_0^{1/(c_1||v_s||_\infty)}c_1^2||v_s||_\infty^{4}t^2dt\Bigg)
^{1/2}\\
&=&\sqrt{\frac{2}{3c_1}}||v_s||_\infty^{1/2}
\end{eqnarray*}
from which we deduce using the second inequality in (\ref{ee1})
$$||v_s||_\infty \leq \frac{3c_1}{2} \cdot c_0^2.$$
This proves the lemma. \hfill $\square$
\\ \\
\textbf{Proof of Theorem~\ref{comhop} continued: }It follows from the
previous lemma that for gradient flow lines of
$\nabla_{g_r}\mathcal{A}_c$ the first factor $v$ remains in the
compact 1-ball around $0$ in the complex plane. Compactness modulo
breaking can now be deduced from the results in 
\cite{cieliebak-gaio-mundet-salamon}. 
However note, that their
arguments simplify in our case. Since our gauge group is abelian
we only need an easy version of Uhlenbeck's compactness theorem.
Moreover, the bubbling analysis can be avoided for the standard
symplectic structure on $\mathbb{C}$ by using the elliptic estimate
$$||v||_{W^{1,2}([-N,N]\times S^1)} \leq c_N
\big(||\bar{\partial}v||_{L^2([-N-1,N+1]\times S^1)}
+||v||_{L^2([-N-1,N+1]\times S^1)}\big)$$
for every $N \in \mathbb{N}$ and a constant $c_N>0$.
\hfill $\square$

\subsection[The maps of Jaffe and Taubes]{The maps of
Jaffe and Taubes}

In \cite{jaffe-taubes} Jaffe and Taubes defined a map
from the moduli space of $N$-vortices on the complex plane to
the $N$-fold symmetric product of the complex plane and
showed that it is bijective. In this subsection we define the analogon of
their map for the gradient flow lines of
$\nabla_{g_r}\mathcal{A}_c$ for all $r \in [0,1]$. We prove
that for $r=1$ the map is bijective. This proves Theorem A
in the introduction. As a biproduct we will obtain the proof of
the compactness statement in Theorem~\ref{comhop}.  

Denote by $\mathfrak{V}^N_r$ the moduli space of $N$-vortices with
respect to the metric $g_r$. Then $\mathfrak{V}^N_1=\mathfrak{V}^N$
the moduli space introduced before. It is useful to write the
map of Jaffe and Taubes
from $\mathfrak{V}^N_r$ to the $N$-fold symmetric product of
the cylinder as the composite of two maps. 
Denote by $\mathfrak{W}^N_r$ the space of distributions
$w$ on the cylinder $\mathcal{Z}$ for which there exists
$N$ not necessarily distinct points $z_j \in \mathcal{Z}$
such that $w$ is smooth outside of $\bigcup_{j=1}^N\{z_j\}$ on
$\mathcal{Z}$ and satisfies the following integro Kazdan-Warner type
problem with singularities and prescribed asymptotic behaviour
\begin{eqnarray}\label{kazdanwarner}
-\Delta w+r^2 e^w+(1-r^2)\int_0^1 e^w dt-1&=&-4 \pi 
\sum_{j=1}^N \delta(z-z_j)\\ \nonumber
\lim_{s\to \pm \infty} w(s,t)&=&0
\end{eqnarray}
where the limit is uniform with respect to the $t$-variable. 
Define the map $\mathfrak{T}^N_r \colon \mathfrak{V}^N_r 
\to \mathfrak{W}^N_r$ by
$$\mathfrak{T}^N_r(v,\eta)=\ln|v|^2$$
and the map $\mathfrak{J}^N_r \colon \mathfrak{W}^N_r \to 
S^N \mathcal{Z}$ by
$$\mathfrak{J}^N_r(w)=[z_1,\ldots,z_N].$$
Note that the composition $\mathfrak{J}^N_r \circ \mathfrak{T}^N_r$
maps a pair $(v,\eta)$ to the zeros of $v$ counted with multiplicity. 
For simplicity of notation we will often drop the index $1$, i.e.
$\mathfrak{J}^N$ means $\mathfrak{J}^N_1$, etc.
We prove the following two theorems.
\begin{thm}\label{taubes}
For every $r \in [0,1]$ and every
$N \in \mathbb{N}$ the map $\mathfrak{T}^N_r$
is bijective.
\end{thm}

\begin{thm}\label{jaffe} 
For every $N \in \mathbb{N}$ the map $\mathfrak{J}^N=\mathfrak{J}^N_1$
is bijective.
\end{thm}
As an easy corollary of the above two theorems we get Theorem A from
the introduction.
\\ \\
\textbf{Proof of Theorem A: } By Theorem~\ref{taubes} and
Theorem~\ref{jaffe} the map $\mathfrak{J}^N \circ \mathfrak{T}^N$
gives a bijection for every $N \in \mathbb{N}$ between
the moduli space of $N$-vortices on the cylinder and the
$N$-fold symmetric product of the cylinder. \hfill $\square$

\begin{rem}
For $r=1$ the problem (\ref{kazdanwarner}) simplifies to 
the following Kazdan-Warner type problem with singularities
(see \cite{kazdan-warner})
\begin{eqnarray}\label{kazdanwarner2}
-\Delta w+e^w-1&=&-4 \pi 
\sum_{j=1}^N \delta(z-z_j)\\ \nonumber
\lim_{s\to \pm \infty} w(s,t)&=&0.
\end{eqnarray}
The bijectivity in Theorem~\ref{jaffe} means that the above problem
has a unique solution. The hard part is to prove existence of a solution.
Our existence proof is based on finite dimensional approximation.
\end{rem}
\textbf{Proof of Theorem~\ref{taubes}: }
Note that for $r \in [0,1]$ the action functional
$\mathcal{A}^r \colon \mathscr{L} \to \mathbb{R}$ defined by
$$\mathcal{A}^r(v,\eta)=\int_0^1 \lambda(v)\partial_t v
+\int_0^1 \langle r \mu(v(t))+(1-r)\bar{\mu}(v),\eta(t)\rangle dt$$
is invariant under the $r$-action of $\mathcal{H}_0$ on
$\mathscr{L}$ and
$$\mathcal{A}^r|_{\mathscr{L}_c}=\mathcal{A}_c.$$
It follows that the gradient flow of $\mathcal{A}^r$ with respect to
the $L^2$-metric on $\mathscr{L}$ are in natural one-to-one correspondence
with gradient flow lines of $\mathcal{A}_c$ with respect to the
$g_r$-metric on $\mathscr{L}_c$ by projection on the Coulomb section.
Note that projection on the Coulomb section does not change the value
of $\ln|v|^2$ and hence we are left with showing the equivalence of
flow line of $\nabla_{g_{L^2}}\mathcal{A}^r$ and solutions of
the problem (\ref{kazdanwarner}). Using the notation
$\bar{\eta}=\int_0^1\eta(t)dt$ for $\eta \in C^\infty(S^1,i\mathbb{R})$
gradient flow lines $(v,\eta)\in C^\infty(\mathcal{Z},\mathbb{C}\times
i\mathbb{R})$
of $\mathcal{A}^r$ with respect to the 
$L^2$-metric solve
\begin{eqnarray}\label{gradr}
\partial_sv+i\partial_t v+ir\eta v+i(1-r)\bar{\eta}v=0\\ \nonumber
\partial_s \eta+r\mu(v)+(1-r)\bar{\mu}(v)=0. 
\end{eqnarray}
It is now an easy exercise to show that 
$\mathfrak{T}^N_r$ is well-defined, i.e. 
$\ln|v|^2$ of solutions of (\ref{gradr}) are solutions of the
integro type Kazdan-Warner type problem with singularities, and
that for each solution of (\ref{kazdanwarner}) there exists a unique
gauge equivalence class of $N$-vortices satisfying (\ref{gradr}),
for details see \cite{jaffe-taubes}. \hfill $\square$
\\ \\
We finally embark on the prove of Theorem~\ref{jaffe}. We first prove
a lemma.

\begin{lemma}\label{inoc}
For every $N \in \mathbb{N}$ the map $\mathfrak{J}^N$ is injective and
its image is open and closed in $S^N \mathcal{Z}$.
\end{lemma}
\textbf{Proof: }
Assume that $w$ and $w'$ are two solution of the problem (\ref{kazdanwarner2})
for the same $N$-tuple of singularities $[z_1,\ldots,z_N]$.
Then its difference $w-w'$ is asymptotically zero and
$\Delta(w-w')\geq e^{w'}(w-w')$. Hence $w=w'$.
\\
The map $\mathfrak{J}^N \circ \mathfrak{T}^N$ is a continuous, one-to-one
map between manifolds of the same dimension. Hence it is open by the
Invariance of domain theorem, see for example 
\cite[Corollary 18.9]{greenberg-harper}.
\\
To show that the image of 
$\mathfrak{J}^N$ is closed, assume
that $w_{\nu} \in \mathfrak{W}^N$ is a sequence such that 
$z_\nu=\mathfrak{J}^N(w_\nu)$ converges to $z \in S^N \mathcal{Z}$
as $\nu$ goes to infinity. For the sequence of flow lines 
$v_\nu=(\mathfrak{T}^N)^{-1}(w_\nu)$ there exists a subsequence
$v_{\nu_j}$ which converges to a broken flow line. But since
$z_\nu$ converges, the limit broken flow line is actually unbroken.
Hence $v_{\nu_j}$ converges to $v \in \mathfrak{V}^N$ and
$$z=\mathfrak{J}^N(\mathfrak{T}^N(v)).$$
Hence the image of $\mathfrak{J}^N$ is closed. \hfill $\square$
\\ \\
The main work lies in the following existence statement for
$1$-vortices. 

\begin{thm}\label{existence}
$\mathfrak{V}^1_r$ is not empty for every $r \in [0,1]$. 
\end{thm}
\textbf{Proof: } Consider the following cylinder action on the
gradient flow lines
$$(v,\eta)(s,t) \mapsto (v,\eta)(s+\sigma,t+\tau),\quad
(\sigma,\tau) \in \mathcal{Z}.$$
$1$-vortices cannot break and hence 
$\mathfrak{V}^1_r/\mathcal{Z}$ is compact by Theorem~\ref{comhop}
for every $r \in [0,1]$. By Proposition~\ref{finiteapp} 
we can identify the flow lines of
$\nabla_{g_0}\mathcal{A}_c$ describing $1$-vortices
with Morse-flow lines on
the finite dimensional approximation $L^1_0$. 
The finite dimensional approximation $L^1_0$ can be identified with
$\mathbb{C}^2$ via its Fourierbasis and the $T^2$-action is given by
$$(e^{i\theta_1},e^{i\theta_2})(z_1,z_2) \mapsto
(e^{i\theta_1}z_1,e^{i(\theta_1+\theta_2)}z_2).$$
Denote by $\mu_{L_0^1}$ the moment map of the circle action
of the first factor in $S^1 \times S^1=T^2$ given by
$$\mu_{L_0^1}(z_1,z_2)=-\frac{i}{2}\big(|z_1|^2+|z_2|^2\big)+\frac{i}{2}.$$
The restriction of the action functional
$$A=\mathcal{A}_c|_{L^1_0 \times i\mathbb{R}}
=\mathcal{A}|_{L^1_0 \times i\mathbb{R} }
\in C^\infty(L^1_0 \times i\mathbb{R})$$
is given by
$$A(v,\eta)=\mathcal{A}_{fl}(v)+\langle \mu_{L_0^1}(v),\eta\rangle.$$
The moduli space $\mathfrak{V}^1_0/\mathcal{Z}$
can now be identied by Proposition~\ref{finiteapp}
with the space of finite energy Morse flow lines of $A$ modulo
$T^2 \times \mathbb{R}$ where the group $\mathbb{R}$ acts
by reparametrisation of flow lines. \\
The space $L^1_0 \times i\mathbb{R}$ is finite dimensional but
still noncompact. Using the results of Appendix~\ref{morselm}
we can homotop our Morse flow lines further to
Morse flow lines on a compact manifold.  
In order to do that note that the function $A$ on $L_0^1 \times i\mathbb{R}$ is
the Lagrange multiplier functional of 
$$H=\mathcal{A}_{fl}|_{\mu_{L_0^1}^{-1}(0)} \in 
C^\infty(\mu^{-1}_{L_0^1}(0)).$$
Hence finite energy Morse flow lines of $A$ can 
be homotoped inside a compact subset of $L^1_0 \times i\mathbb{R}$
to Morse flow lines of $H$.
The manifold $\mu^{-1}_{L_0^1}(0)$ is the three sphere $S^3$,
the circle action of the first factor in $T^2=S^1\times S^1$
is the Hopf fibration $S^3 \to S^2$ and the circle action of
the second factor in $T^2$ acts by rotation on the two-sphere
$S^2$. The function $H$ induces on $S^2$ the height function. 
In particular, the action of $T^2 \times \mathbb{R}$ on
the Morse flow lines of $H$ is free and the quotient consists
of exactly one point. 
\\
The upshot of our construction is that we can homotop the moduli
space $\mathfrak{V}^1_r/\mathcal{Z}$ by a compact homotopy
to a point. During this homotopy the Fredholm index is unchanged
by Proposition~\ref{lagmorbo}.
We are now in position to show that $\mathfrak{V}^1_r$ is nonempty
for every $r \in [0,1]$. Assume the contrary. Then we apply the
abstract perturbation theory of \cite{fukaya-ono,li-tian,liu-tian,ruan,siebert}
to our compact homotopy. Actually, since we do not have to compactify
our moduli spaces by broken flow lines containing bubble trees, 
the more elementary theory of \cite{cieliebak-mundet-salamon}
is already sufficient. What we obtain is a compact branched manifold
containing just one boundary point of weight one. But
such an object does not exist. Hence
$\mathfrak{V}^1_r$ is nonempty for every $r \in [0,1]$. \hfill $\square$
\\ \\
\textbf{Proof of Theorem~\ref{jaffe}: }
It follows from Floer's gluing construction and
Theorem~\ref{existence} that $\mathfrak{V}^N$ is not empty for 
every $N \in \mathbb{N}$. Hence 
$\mathfrak{J}^N(\mathfrak{W}^N)=\mathfrak{J}^N\circ
\mathfrak{T}^N(\mathfrak{V}^N)$ is not empty in 
$S^N \mathcal{Z}$. Since the image im $\mathfrak{J}^N$ is open and
closed by Lemma~\ref{inoc} and $S^N\mathcal{Z}$ is connected it follows that
$\mathfrak{J}^N$ is surjective. Since again by Lemma~\ref{inoc}
$\mathfrak{J}^N$ is injective
the theorem follows. \hfill $\square$

\section[Further directions]{Further directions}\label{Further}
\subsection[The symplectic vortex equations and Givental's toric 
map spaces]{The symplectic vortex equations and Givental's toric
map spaces}
Instead of the circle action on $\mathbb{C}$ we can study
more generally linear torus actions on a complex vector space.
Assume that for $k\leq n$ the torus 
$T^k=\{e^{iv}: v \in \mathbb{R}^k\}$ 
acts on the complex vector space
$\mathbb{C}^n$ via the action
$$\rho(e^{iv})z=e^{iAv}z, \quad z \in \mathbb{C}^n, \,\,v \in \mathbb{R}^k$$
for some $(n\times k)$-matrix $A$ with 
integer entries. We endow the 
Lie algebra of the torus
$$\mathrm{Lie}(T^k)=\mathfrak{t}^k=i\mathbb{R}^k$$
with its standard inner product. 
The action of the torus
on $\mathbb{C}^n$ is Hamiltonian with respect to the standard symplectic
structure $\omega=\sum_{i=1}^n dx_i \wedge dy_i$. 
Denoting by $A^T$ the transposed matrix of $A$ a moment map 
$\mu \colon \mathbb{C}^n \to \mathfrak{t}^k$ is given by
\begin{equation}\label{moment1}
\mu(z)=-i A^Tw, \quad w=\frac{1}{2}\left(\begin{array}{c}
|z_1|^2 \\
\vdots\\
|z_n|^2
\end{array}\right),
\end{equation}
i.e.
$$d\langle \mu, \xi \rangle=\iota_{X_\xi}\omega, \quad \xi \in 
\mathfrak{t}^k$$
for the vector field $X_\xi$ on $\mathbb{C}^n$ given by the
infinitesimal action
$$X_\xi(z)=\dot{\rho}(\xi)(z), \quad z \in \mathbb{C}^n.$$
We assume the following hypothesis, 
\begin{description}
 \item[(H)] The moment map $\mu$ is proper and $T^k$ acts freely on $\mu^{-1}(\tau)$.
\end{description}
It follows from (H) that the Marsden-Weinstein quotient
$$\mathbb{C}^n//T^k=\mu^{-1}(\tau)/T^k$$
is a compact symplectic manifold of dimension
$$\mathrm{dim}(\mathbb{C}^n//T^k)=2(n-k),$$
where the symplectic structure is
induced from the standard symplectic structure on $\mathbb{C}^n$.
\\
Let $\mathscr{L}$ be the loop space
$$\mathscr{L}:=C^\infty(S^1,\mathbb{C}^n \times \mathfrak{t}^k).$$
The gauge group
$$\mathcal{H}=C^\infty(S^1,T^k)$$
acts on $\mathscr{L}$ by
$$h_*(v,\eta)=(\rho(h)v,\eta-h^{-1}\partial_t h), \quad
h \in \mathcal{H},\,\,(v,\eta) \in \mathscr{L}.$$
Recall Floer's action functional 
$\mathcal{A}_{fl}\colon C^\infty(S^1,\mathbb{C}^n) \to \mathbb{R}$
given by
$$\mathcal{A}_{fl}(v)=\int_0^1 \lambda(v)(\partial_t v)$$
where $\lambda$ denotes the Liouville 1-form 
$$\lambda=\sum_{i=1}^n y_i dx_i, \quad d\lambda=-\omega.$$
The Moment action functional 
$\mathcal{A} \colon \mathscr{L} \to \mathbb{R}$ 
is defined by
$$\mathcal{A}(v,\eta):=\mathcal{A}_{fl}(v)+
\int_0^1\langle \mu(v(t))-\tau,\eta(t)
\rangle dt.$$
Again one may think of $\eta$ in the
second integral as a
Lagrange multiplier. In particular,
the critical points of $\mathcal{A}$ are the critical points
of Floer's action on the constraint $\mu^{-1}(\tau)$.
\\
The gradient flow lines
of $\mathcal{A}$ with respect to the $L^2$-metric $g_{L^2}$ on
$\mathscr{L}$ are solutions $(v,\eta)\in C^\infty(\mathbb{R} \times S^1,
\mathbb{C}^n\times \mathfrak{t}^k)$ of
\begin{eqnarray*}
\partial_s v+i\partial_t v+i\dot{\rho}(\eta)v=0\\
\partial_s \eta+\mu(v)=\tau.
\end{eqnarray*}
These are examples of the symplectic vortex equations on the cylinder in 
temporal gauge.

In this setting the symplectic vortex equations can be homotoped again
via a Floer homotopy compact up to breaking of flow lines to
Morse problems on  finite dimensional compact manifolds.
If one considers a finite dimensional Fourierapproximation
$L$ of the loop space $C^\infty(S^1, \mathbb{C}^n)$ then
the $T^k$-action on $\mathbb{C}^n$ induces a $T^k$-action
on $L$ by coefficientwise multiplication. This action is again
Hamiltonian with moment map $\mu_L$ normalized such that
$\mu_L(0)=0$. The finite dimensional compact manifolds
we end up with are 
$$G_L=\mu_L^{-1}(\tau)/T^k$$
and the Morse function is again Floer's action functional
restricted to $G_L$. The spaces $G_L$ are known as
Givental's toric map spaces. They were
introduced by Givental in \cite{givental2} and
studied by different authors in \cite{givental,iritani,vlassopoulos}.     

\subsection[Warped product metrics and Chern-Simons Vortices]{
Warped product metrics and Chern-Simons Vortices}
 
If we consider the same action functional 
$\mathcal{A} \colon \mathscr{L} \to \mathbb{R}$
as for the vortex equations but instead of the standard
flat $L^2$-metric a warped product metric $g_\gamma$
for a smooth function $\gamma\colon [0,\infty) \to (0,\infty)$
given by
$$g_\gamma(v,\eta)\big((\hat{v}_1,\hat{\eta}_1),
(\hat{v}_2,\hat{\eta}_2)\big)=
\int_0^1\langle \hat{v}_1,\hat{v}_2 \rangle dt
+\int_0^1 \gamma(|v|)^2\langle \hat{\eta}_1,
\hat{\eta}_2 \rangle dt$$
for $(v,\eta) \in \mathscr{L}$ and $(\hat{v}_1,\hat{\eta}_1),
(\hat{v}_2,\hat{\eta}_2) \in T_{(v,\eta)}\mathscr{L}$
we obtain the following gradient equations
for 
$(v,\eta) \in C^\infty(\mathbb{R}\times S^1, \mathbb{C}\times i\mathbb{R})$
\begin{eqnarray*}
\partial_s v+i\partial_t v+i\eta v=0\\
\partial_s \eta+\frac{1}{\gamma(|v|)^2} \mu(v)=0.
\end{eqnarray*}
In particular, if we choose
$$\gamma(r)=\frac{1}{r}$$
we obtain
\begin{eqnarray*}
\partial_s v+i\partial_t v+i\eta v=0\\
\partial_s \eta+|v|^2\mu(v)=0.
\end{eqnarray*}
These are the selfduality equations for
the Chern-Simons vortices on the cylinder discovered
by Hong-Kim-Pac and Jackiw-Weinberg, see
\cite{hong-kim-pac, jackiw-weinberg}. We refer the reader
to the excellent textbook of Y.\,Yang \cite{yang} for a detailed
treatment of this equation. This textbook may also serve as a guide
to the corresponding literature. 

Note that for this choice of $\gamma$ the metric $\gamma_v$ becomes
singular if $v$ goes to zero. One may think of this as a continuum
of ``critical points at infinity'' for each $(0,\eta)$
where $\eta$ is a smooth loop in the Lie algebra $i\mathbb{R}$.
In particular, the action functional $\mathcal{A}$ takes on
the set of ``critical points at infinity'' every value in 
$\mathbb{R}$. The energy of flow lines which converge at
one end to a ``critical point at infinity'' can therefore be
any value in $\mathbb{R}$, in contrast to the classical vortex
equations where the energy of finite energy flow lines was
quantized. Such solutions are called in the physics literature
``nontopological solutions''. We again refer to the textbook
of Y.\,Yang \cite{yang} and the literature cited therein 
for a detailed treatment of nontopological solutions.
For compact Riemann surfaces existence of ``nontopological solutions''
was proved by Tarantello and Ding-Jost-Li-Peng-Wang
\cite{ding-jost-li-peng-wang,tarantello}. In this case the
``nontopological solutions'' are characterized by the property
that $v$ converges to $0$ under the adiabatic limit
obtained by letting the Chern-Simons coupling parameter tend to zero.  

On the finite dimensional
approximations the Chern-Simons vortices are flow lines of
a Lagrange multiplier functional with respect to a warped product
metric. ``Critical points at infinity'' are responsible
for the failure of the Palais-Smale condition discussed in the
appendix. So the study of flow lines of Lagrange multiplier functionals with
respect to a warped product metric is a finite dimensional analogon
of the Abelian Chern-Simons-Higgs theory and should lead to
a deeper understanding of the phenomenons occuring in this theory.

\appendix \label{ML}

\section[Morse functions with Lagrange multipliers]{Morse functions with
Lagrange multipliers}\label{morselm}

Assume that $M$ is a finite dimensional
manifold and $V$ is a finite dimensional real vector space.
It is well known from basic calculus that critical points
of a smooth function $f \in C^\infty(M)$ satisfying a constraint
given by the zero set of a smooth function $h \in C^\infty(M,V)$
can be found by considering the Lagrange multiplier functional
$F \in C^\infty(M \times V^*)$, where $V^*$ is the dual vector space of
$V$, given by
$$F(x,v^*)=f(x)+v^*(h(x)).$$
If $0$ is a regular value of $h$ then there is a natural one-to-one
correspondence between critical points of $F$ and critical points
of $f|_{h^{-1}(0)}$. However, Morse flow lines of $F$ and
Morse flow lines of $f|_{h^{-1}(0)}$ may be quite different. 
Even if $h^{-1}(0)$ is compact it is not a priori clear that
the moduli spaces of flow lines of $F$ are compact modulo breaking
since $F$ is neither bounded from above nor below. However, we will
show that if $h$ is locally proper around $0$, then
$F$ satisfies the Palais-Smale condition from which we can deduce
that flow lines of $F$ remain in a compact subset of the
noncompact manifold $M \times V^*$.

A first possibility to homotop Morse flow lines of $F$ to
Morse flow lines of $f|_{h^{-1}(0)}$ would be the adiabatic limit method.
For a fixed Riemannian metric $g_M$ on $M$ and a fixed Riemannian
metric $g_{V^*}$ on $V^*$, induced from a Euclidean scalar product on $V$
we consider the family of metrics $g_\epsilon$ on
$M \times V^*$ for $\epsilon \in (0,1]$
$$g_\epsilon=g_M \oplus \epsilon^2 g_{V^*}.$$
If $\epsilon$ goes to zero, then the gradient flow lines of $F$
with respect to the metric $g_\epsilon$ converge to gradient flow
lines of $f|_{h^{-1}(0)}$ with respect to the metric
$g_M|_{h^{-1}(0)}$. If a generalized implicit function theorem
as in \cite{dostoglou-salamon} can be established at $\epsilon=0$,
then this would lead to a homotopy compact modulo breaking
between the two moduli spaces of gradient flow lines. 

In this section we will pursue another approach. We will consider
a homotopy of $f$ and the Riemannian metric $g_M$ 
such that $f|_{h^{-1}(0)}$ is unchanged during the
homotopy but the normal derivatives of $\nabla_{g_M}f$ at 
$h^{-1}(0)$ are homotoped
to zero. Since $f|_{h^{-1}(0)}$ is fixed the critical points
of the Lagrange multiplier functional can be canonically identified
with the set of critical points of $f|_{h^{-1}(0)}$ during the
whole homotopy. If the normal derivatives of $\nabla_{g_M} f$
at $h^{-1}(0)$ vanish then the moduli space of Morse flow lines
of $f|_{h^{-1}(0)}$ is canonically contained in the moduli space
of flow lines of $F$. We will prove that for special choices
of $f$ and $g_M$ there are no other flow lines of $F$. The main
idea is to choose $g_M$ in such a way that a tubular neighbourhood
of $h^{-1}(0)$ in $M$ becomes very huge and then prove that finite
energy flow lines have to remain in this tubular neighbourhood.
\\

It is natural to formulate our main theorem in the language of
Morse-Bott functions. In order to fix notation we 
recall briefly its definition.
A function $F$ on a finite dimensional manifold $M$ 
is called Morse-Bott if
the critical set
is a submanifold of $M$ and for each $x \in \mathrm{crit}(F)$
we have
$$T_x \mathrm{crit}(F)=\mathrm{ker}H_F(x)$$
where $H_F(x)$ is the Hessian of $F$ at $x$. 
It is well known that the Morse-Bott condition implies that 
flow lines which remain in a compact set of $M$
converge at both ends exponentially fast to
critical points of $F$. For a Riemannian metric $g$ on $M$ we denote by
$\mathcal{M}(F,g)$
the moduli space of finite energy flow lines of $\nabla_g F$.
\\ \\
The main theorem of this section can now be statet in the following way.
 
\begin{thm}\label{lagmor}
Let $M$ be a finite dimensional manifold, $\Gamma$ be a Lie group acting
on $M$,
and let $(V,\langle\,,\,\rangle)$ be a finite dimensional 
Euclidean vector space.
Assume that $g_M$ is a $\Gamma$-invariant, geodesically complete Riemannian
metric on $M$,
$h \in C^\infty(M,V)$
and $f \in C^\infty(M)$ are $\Gamma$-invariant
functions satisfying the following conditions.
\begin{itemize}
 \item $0$ is a regular value of $h$,
 \item $h$ is locally proper around $0$, i.e. 
  there exists an open neighbourhood $V_0$ of $0$ in $V$ 
  such that $h^{-1}(\mathrm{cl}(V_0))$ is compact,
 \item the restriction of $f$ to 
  the compact manifold $h^{-1}(0)$ is Morse-Bott. 
\end{itemize}
Denote by $V^*$ the dual vector space of $V$ and let
$\Gamma$ act on $M \times V^*$ by $\gamma(x,v^*)=(\gamma x, v^*)$
for $\gamma \in \Gamma$ and $(x,v^*) \in M \times V^*$.
Then there exists a smooth family of 
$\Gamma$-invariant Morse-Bott functions
$F_r \in C^\infty(M \times V^*)$ for $r \in [0,1]$ and a smooth 
family of $\Gamma$-invariant 
Riemannian metrics $g_r$ on $M \times V^*$
satisfying
$$F_0(x,v^*)=f(x)+v^*(h(x)),
  \quad g_0=g_M \oplus g_{V^*}$$
where $g_{V^*}$ is the metric on $V^*$ induced from
the scalar product $\langle \,,\,\rangle$ on $V$, 
such that the following conditions are satisfied. 
\begin{description}
 \item[(i)] The inclusion $\iota \colon h^{-1}(0) \to M \times V^*,
 x \mapsto (x,0)$ induces a bijection
 $$\iota_* \colon \mathcal{M}(f|_{h^{-1}(0)},g_M|_{h^{-1}(0)})
  \to \mathcal{M}(F_1,g_1)$$
 defined by
 $$\iota_*y(s)=\iota(y(s)), \quad 
  y \in \mathcal{M}(f|_{h^{-1}(0)},g_M|_{h^{-1}(0)}),\,\,s \in \mathbb{R}.$$
\item[(ii)] For $r \in [0,1]$ there exists
  a smooth family of diffeomorphism  $\phi_r \colon \mathrm{crit} F_0
  \to \mathrm{crit}F_r \subset M \times V^*$.
 \item[(iii)] There exists a compact set $K \subset M \times V^*$ such 
  that 
  $$\Bigg\{y(\sigma):y \in \bigcup_{r \in [0,1]} 
  \mathcal{M}(F_r,g_r),\,\,\sigma \in \mathbb{R}\Bigg\}
  \subset K.$$
 \end{description}
\end{thm}
\textbf{Proof:} We prove the theorem in seven steps.
\\ \\
\textbf{Step 1 (Neighbourhood of the constraint): }
\emph{There exists an open neighbourhood
$V_1$ of $0$ in $V$, a $\Gamma$-invariant open neighbourhood
$U$ of $h^{-1}(0)$ in $M$, and a $\Gamma$-equivariant diffeomorphism
$$\phi\colon h^{-1}(0)\times V_1 \to U$$
where $\Gamma$ acts on $h^{-1}(0)\times V_1$ by
$\gamma(x,v)=(\gamma x, v)$ for $\gamma \in \Gamma$ and
$(x,v) \in h^{-1}(0)\times V_1$ such that
\begin{equation}\label{tubu}
h(\phi(x,v))=v, \quad (x,v) \in h^{-1}(0) \times V_1.
\end{equation}}
\\
Since $0$ is a regular value of $h$ there exists an open neighbourhood $U_0$
of $h^{-1}(0)$ such that $dh(y)$ is surjectiv for every $y \in U_0$.
For $v \in V$ define the vector field $\xi_v$ on $U_0$
by the conditions
$$dh(y)\xi_v(y)=v, \quad \xi_v(y) \in \mathrm{ker} dh(y)^\perp,
\quad y \in U_0$$ 
where $\mathrm{ker}dh(y)^\perp$ denotes the orthogonal complement
of the kernel of $dh(y)$ with respect to the metric $g_M$.
Since $h^{-1}(0)$ is compact there exists an open neighbourhood
$V_1$ of $0$ in $V$ such that for each $v \in V_1$ and for
each $x \in h^{-1}(0)$ there exists a unique solution 
$y_{x,v} \in C^\infty([0,1],U)$ of the problem
$$y_{x,v}(0)=x,\quad 
\partial_t y_{x,v}(t)=\xi_v(y_{x,v}(t)),\quad t \in [0,1],$$
and the map
$$\phi(x,v)=y_{x,v}(1), \quad (x,v) \in h^{-1}(0) \times V_1$$
is a diffeomorphism. Set $U=\phi(h^{-1}(0) \times V_1)$.
Since $h$ and $g_M$ are $\Gamma$-invariant it follows that
$\gamma (y_{x,v})=y_{\gamma x,v}$ for $\gamma \in \Gamma$ and
$(x,v)\in h^{-1}(0) \times V_1$. Hence $U$ is $\Gamma$-invariant
and $\phi \colon h^{-1}(0)\times V_1 \to U$ is a
$\Gamma$-equivariant diffeomorphism.
Moreover, we compute
\begin{eqnarray*}
h(\phi(x,v))&=&h(y_{x,v}(1))\\
&=&h(y_{x,v}(0))+\int_0^1 \frac{d}{dt}h(y_{x,v}(t))dt\\
&=&h(x)+\int_0^1 dh(y_{x,v}(t))\partial_t y_{x,v}(t)dt\\
&=&\int_0^1 dh(y_{x,v}(t))\xi_v(y_{x,v}(t))dt\\
&=&\int_0^1 v dt\\
&=&v.
\end{eqnarray*}
This proves (\ref{tubu}) and hence Step 1.
\\ \\
\textbf{Step 2: } \emph{We construct the homotopies.}
\\ \\
In this step we construct a $\Gamma$-invariant function
$f_1 \in C^\infty(M)$ and a $\Gamma$-invariant metric $g_{M,1}$ on $M$.
We then set for $r \in [0,1]$
$$f_r=(1-r)f+rf_1,\quad g_{M,r}=(1-r)g_M+rg_{M,1}$$
and define the homotopy of functions $F_r \in C^\infty(M\times V^*)$
by
$$F_r(x,v^*)=f_r(x)+v^*(h(x)), \quad (x,v^*) \in M\times V^*$$
and the homotopy of metrics $g_r$ on $M \times V^*$ by
$$g_r=g_{M,r}\oplus g_{V^*}.$$

Choose a small number $\delta>0$ such that the open $\delta$-ball
$B_\delta=\{v \in V:||v|| <\delta\}$ is contained in the neighbourhood
$V_1$ of $0$ in $V$ constructed in Step 1. Moreover, since $h$ is
locally proper at $0$, we may assume that
\begin{equation}\label{sorry}
h(x)>\delta \Rightarrow x \notin \phi(h^{-1}(0)\times B_\delta).
\end{equation}
Choose further a
cutoff function $\hat{\beta} \in C^\infty([0,\delta),[0,1])$
such that $\hat{\beta}|_{[0,\delta/2]}=1$ and 
$\hat{\beta}|_{[3\delta/4,\delta)}=0$.
Denote by $\pi_1\colon h^{-1}(0) \times B_\delta \to h^{-1}(0)$ and
by $\pi_2 \colon h^{-1}(0)\times B_\delta \to B_\delta$ the projection
to the first, respectively the second, factor. We will use the
following $\Gamma$-invariant cutoff function on $M$ given by
$$\beta(x)=\left\{\begin{array}{cc}
\beta\big(|\pi_2(\phi^{-1}(x))|\big)
& x \in \phi(h^{-1}(0)\times B_\delta)
\\
0 & x \notin \phi(h^{-1}(0)\times B_\delta).
\end{array}\right.$$
Define the function
$f_1 \in C^\infty(M)$ by
$$f_1(x)=\left\{\begin{array}{cc}
\beta(x)f\big(\pi_1(\phi^{-1}(x))\big)+
\big(1-\beta(x)\big)f(x)
& x \in \phi(h^{-1}(0)\times B_\delta)
\\
f(x) & x \notin \phi(h^{-1}(0)\times B_\delta).
\end{array}\right.$$
Set
\begin{equation}\label{energybound}
C:=\max_{x \in h^{-1}(0)}\{f(x)\}-\min_{x \in h^{-1}(0)}\{f(x)\}
\end{equation}
and choose a constant
\begin{equation}\label{distance}
\kappa>\frac{16C}{\delta^2}.
\end{equation}
Let $g_\kappa$ be the product metric on $h^{-1}(0)\times B_\delta$
$$g_\kappa=g_M|_{h^{-1}(0)} \oplus \kappa^2 g_{B_\delta}$$
where $g_{B_\delta}$ is the standard euclidean metric on the ball
$B_\delta \subset V$.
We are now able to define the metric $g_{M,1}$ on $M$ by the
formula
$$g_{M,1}(x)=\left\{\begin{array}{cc}
\beta(x)(\phi_*g_\kappa)(x)+
\big(1-\beta(x)\big)g_M(x)
& x \in \phi(h^{-1}(0)\times B_\delta)
\\
g_M(x) & x \notin \phi(h^{-1}(0)\times B_\delta).
\end{array}\right.$$
\textbf{Step 3: }\emph{The trace of each finite
energy flow line $y \in C^\infty(\mathbb{R},M\times V^*)$
of $\nabla_{g_1} F_1$ is contained in
$\phi(h^{-1}(0) \times B_{\delta/2})\times V^*$.}
\\ \\
First note that if $x \in M \setminus \phi(h^{-1}(0)
\times B_{\delta/4})$ and $v^* \in V^*$ then it follows
from (\ref{tubu}) and (\ref{sorry}) that
\begin{equation}\label{estad1}
||\nabla_{g_1} F_1(x,v^*)||\geq ||h(x)||
\geq \frac{\delta}{4}.
\end{equation}
Observe further that the energy of a finite energy flow line is bounded from
above by the constant $C$ introduced in (\ref{energybound}), i.e.
\begin{equation}\label{estad2}
\int_{-\infty}^\infty||\nabla_{g_1} F_1(y(s))||_{g_1}^2 ds \leq C.
\end{equation}
Now assume by contradiction that there exists $\sigma \in \mathbb{R}$
such that 
\begin{equation}\label{adestcont}
y(\sigma) \notin \phi(h^{-1}(0) \times B_{\delta/2}) \times V^*.
\end{equation}
Denote by $\tau(\sigma)>\sigma$ the real number
$$\tau(\sigma):=\min \{s \in \mathbb{R}: y(s)
\notin \phi(h^{-1}(0) \times B_{\delta/4}) \times V^*\}.$$ 
Note that $\tau(\sigma)$ is finite, since the energy of the flow line 
$y$ is assumed to be finite and the critical points of 
$F_1$ lie in $h^{-1}(0) \times V^*$.
Denoting by $\mathrm{dist}_{g_1}(\cdot,\cdot)$ 
the distance with respect to the metric
$g_1$ we estimate using (\ref{distance}), (\ref{estad1}), and
(\ref{estad2})
\begin{eqnarray*}
C &<&\frac{\delta^2 \kappa}{16}\\
&\leq&\frac{\delta}{4}\cdot \mathrm{dist}_{g_1} 
\big(y(\sigma),y(\tau(\sigma))\big)\\
&\leq&\frac{\delta}{4}\int_\sigma^{\tau(\sigma)}||\partial_s y(s)||_{g_1}ds\\
&=&\frac{\delta}{4}\int_\sigma^{\tau(\sigma)}||\nabla_{g_1} F_1(y(s))||
ds\\
&\leq&\int_\sigma^{\tau(\sigma)}||\nabla_{g_1} F_1(y(s))||^2 ds\\
&\leq&\int_{-\infty}^\infty||\nabla_{g_1} F_1(y(s))||^2 ds\\
&\leq& C.
\end{eqnarray*}
This contradiction shows that (\ref{adestcont}) cannot hold which proves 
Step 3.
\\ \\
\textbf{Step 4: }\emph{The trace of each finite
energy flow line $y=(x,v^*) \in C^\infty(\mathbb{R},M\times V^*)$
of $\nabla_{g_1} F_1$ is contained in
$h^{-1}(0)\times \{0\}$.}
\\ \\
It follows from Step 3 that  $x$ is contained in the image of $\phi$. 
Denoting
$$(q,w)=\phi^{-1}(x) \in C^\infty(\mathbb{R},h^{-1}(0)\times B_\delta)$$
we observe that the triple $(q,w,v^*)$ is a flow line of the
function
$$F(q,w,v^*)=f(q)+v^*(w)$$
with respect to the metric
$$g_M|_{h^{-1}(0)} \oplus \kappa^2 g_{B_\delta} \oplus g_{V^*}.$$
Denote by $\Lambda \colon V^* \to V$ the isomorphism induced from the
euclidean scalar product on $V$. 
Flow lines of $F$ are solutions of the following ODE
\begin{eqnarray}\label{norlag}
\partial_s q+\nabla_{g_M|_{h^{-1}(0)}}f(q)=0 \\ \nonumber
\partial_s w+\frac{1}{\kappa}\Lambda v^*=0 \\ \nonumber
\partial_s v^*+\Lambda^{-1}w=0.
\end{eqnarray}
It follows from the two last equations in (\ref{norlag}) that
there exist $w_0, w_1 \in V$ such that
$$w(s)=w_0 \exp\bigg(\frac{s}{\sqrt{\kappa}}\bigg)
+w_1 \exp\bigg(-\frac{s}{\sqrt{\kappa}}\bigg).$$
Since the energy of the flow line $y$ is assumed to be finite it follows
that $w_0=w_1=0$ and hence
$$w(s)=0, \quad v^*(s)=0, \quad s\in \mathbb{R}.$$
This proves Step 4.
\\ \\
\textbf{Step 5 (Uniform Palais-Smale condition):}
 \emph{There exists a geodesically complete Riemannian  metric $g_{PS}$ on
  $M \times V^*$, a
  compact set $K_0 \subset M \times V^*$ and a constant $\epsilon>0$
  such that for $y \in (M \times V^*) \setminus K_0$ and
  $r \in [0,1]$
 \begin{equation}\label{PS} 
  \nabla_{gr} F_r(y) \neq 0, \quad
  ||\nabla_{g_r}F_r (y)||^2_{g_r} \geq \epsilon||\nabla_{g_r}F_r(y)||_{g_{PS}}
 \end{equation}
 where $||\,||_{g}$ denotes the norm induced from the metric $g$.}
\\ \\
We choose $g_{PS}=g_M \oplus g_{V^*}=g_0$. 
Then $g_{PS}$ is geodesically complete
by assumption. For $x \in M$ and $r \in [0,1]$
we denote by $dh(x)^{*_r} \colon V^* \to T_xM$ 
the adjoint of $dh(x)$ with respect 
to the inner products $g_{M,r}(x)$ on $T_x M$ and
$\langle\,,\,\rangle$ on $V^*$. With respect to the natural splitting
$T_{(x,v^*)}(M \times V*) \cong T_x M \times V^*$ for 
$(x,v^*) \in M \times V^*$ the gradient of $F_r$ reads
\begin{eqnarray}\label{gradPS}
\nabla_{g_r}F_r(x,v^*)&=&\left(\begin{array}{c}
\nabla_{g_{M,r}}f_r(x)+\nabla_{g_{M,r}}(v^* \circ h)(x)
\\
h(x)
\end{array}\right)\\ \nonumber
&=&
\left(\begin{array}{c}
\nabla_{g_{M,r}}f_r(x)+dh^{*_r}(x)v^*
\\
h(x)
\end{array}\right).
\end{eqnarray}
Since $0$ is a regular value of $h$ and 
$h^{-1}(\mathrm{cl}(V_0))$ is compact we can find
an open neighbourhood $V_0'$ of $0$ in $V$ satisfying
$V_0' \subset V_0$ such that $dh(x)^{*_r}$ is injectiv for
every $x \in \mathrm{cl}(V_0')$ and every $r \in [0,1]$. Set
$$\epsilon':=\min_{v \in V \setminus V_0'}||v||>0.$$
Since the family of injective maps $dh(x)^{*_r}$ depends smoothly
on the compact parameter 
$(x,r) \in h^{-1}(\mathrm{cl}(V_0')) \times [0,1]$ there exists 
a compact subset $W \in V^*$ such that
\begin{equation}\label{estPS}
||\nabla_{g_{M,r}} f_r(x)+dh(x)^{*_r}v^*||_{g_{M,r}} \geq \epsilon',\quad
v^* \in V^* \setminus W,\,\, x \in h^{-1}(\mathrm{cl}(V_0')),\,\, r\in [0,1].
\end{equation}
We set
$$K_0=h^{-1}(\mathrm{cl}(V_0')) \times W.$$
Then $K_0$ is compact and we claim that
\begin{equation}\label{PS1}
 ||\nabla_{g_{M,r}}f_r(x)+dh(x)^{*_r}v^*||_{g_{M,r}}+||h(x)||
\geq \epsilon', \quad y=(x,v^*) \in (M \times V^*)\setminus
K_0.
\end{equation}
To prove the claim we first assume that 
$x \notin h^{-1}(\mathrm{cl}(V_0'))$. We then estimate 
$$||\nabla_{g_{M,r}}f_r(x)+dh(x)^{*_r}v^*||_{g_{M,r}}+||h(x)||
\geq ||h(x)|| \geq \epsilon'$$
by the definition of $\epsilon'$. Now assume that
$x \in h^{-1}(\mathrm{cl}(V_0'))$ but
$v^* \notin W$. We estimate in this case using 
(\ref{estPS})
$$
||\nabla_{g_{M,r}}f_r(x)+dh(x)^{*_r}v^*||_{g_{M,r}}+||h(x)||\geq
||\nabla_{g_r} f_r(x)+dh(x)^{*_r}v^*||_{g_{M,r}} \geq \epsilon'.
$$
This proves (\ref{PS1}).
\\
Using (\ref{gradPS}) and (\ref{PS1}) we estimate for 
$y=(x,v^*) \in (M \times V^*) \setminus K_0$
$$||\nabla_{g_r}F_r(y)||_{g_r}\geq 
\frac{1}{\sqrt{2}}\bigg(
||\nabla_{g_{M,r}}f_r(x)+dh(x)^{*_r}v^*||_{g_{M,r}}+||h(x)||\bigg)
\geq \frac{\epsilon'}{\sqrt{2}}>0$$
which implies the first inequality in (\ref{PS}). To prove the second one
we observe that
since the metrics $g_{M,r}$ differ from $g_M$ only on a compact subset
of $M$ the metrics $g_{M,r}$ are equivalent for every $r \in [0,1]$,
i.e. there exists a constant $c \geq 1$ such that
\begin{equation}\label{metricPS}
\frac{1}{c^2}g_M \leq g_{M,r} \leq c^2 g_M, \quad r \in [0,1].
\end{equation}
Using (\ref{gradPS}), (\ref{PS1}), and (\ref{metricPS}) we estimate
for 
$y=(x,v^*) \in (M \times V^*) \setminus K_0$
\begin{eqnarray*}
||\nabla_{g_r}F_r(y)||^2_{g_r}
&=&||\nabla_{g_{M,r}}f_r(x)+dh^{*_r}(x)v^*||_{g_{M,r}}^2+||h(x)||^2\\
&\geq&\frac{1}{2}\bigg(
||\nabla_{g_{M,r}}f_r(x)+dh^{*_r}(x)v^*||_{g_{M,r}}+||h(x)||\bigg)^2\\
&\geq&\frac{\epsilon'}{2}
\bigg(\frac{1}{c}
||\nabla_{g_{M,r}}f_r(x)+dh^{*_r}(x)v^*||_{g_M}+||h(x)||\bigg)\\
&\geq&\frac{\epsilon'}{2^{3/2}c}||\nabla_{g_r}F_r(y)||_{g_{PS}}.
\end{eqnarray*}
Hence the second inequality in (\ref{PS}) follows with
$\epsilon=\epsilon'/2^{3/2}c$. This proves Step 5.
\\ \\
\textbf{Step 6:} \emph{We prove (iii).}
\\ \\
Let $y \in \bigcup_{r \in [0,1]}\mathcal{M}(F_r,g_r)$. Let
$K_0 \subset M \times V^*$ be the compact set found in Step 5. We estimate
for each $\sigma \in \mathbb{R}$ the distance 
$\mathrm{dist}_{PS}(y(\sigma),K_0)$ between $y(\sigma)$ and $K_0$
with respect to the Palais-Smale metric $g_{PS}$ found in Step 5.
We abbreviate
$$m:=\max_{\substack{x \in K_0,\\r \in [0,1]}}F_r(x)-
\min_{\substack{x \in K_0,\\r \in [0,1]}}F_r(x).
$$
Since the Morse flow line $y$ has finite energy it follows from 
(\ref{PS}) that for each $\sigma \in \mathbb{R}$ the set
$\{s \geq \sigma: y(s) \in K_0\}$ is nonempty. We set
$$\tau(\sigma)=\inf\{s \geq \sigma: y(s) \in K_0\}.$$
Using (\ref{PS}) and the gradient equation we estimate
\begin{eqnarray*}
\mathrm{dist}_{PS}(y(\sigma),K_0)&\leq&
\int_{\sigma}^{\tau(\sigma)}||\partial_s y(s)||_{g_{PS}} ds\\
&=&\int_{\sigma}^{\tau(\sigma)}||\nabla_{g_r} F_r(y(s))||_{g_{PS}} ds\\
&\leq&\frac{1}{\epsilon}
\int_{\sigma}^{\tau(\sigma)}||\nabla_{g_r} F_r(y(s))||^2_{g_r}ds\\
&\leq&\frac{1}{\epsilon}
\int_{-\infty}^{\infty}||\nabla_{g_r} F_r(y(s))||^2_{g_r}ds\\
&=&-\frac{1}{\epsilon}\int_{-\infty}^\infty
g_r(y(s))(\nabla_{g_r}F_r(y(s)),\partial_s y(s))ds\\
&=&-\frac{1}{\epsilon}\int_{-\infty}^\infty dF_r(y(s))\partial_s y(s)ds\\
&=&-\frac{1}{\epsilon}\int_{-\infty}^\infty\frac{d}{ds} F_r(y(s))ds\\
&\leq&\frac{1}{\epsilon}\bigg(\limsup_{s \to -\infty}F_r(y(s))-
\liminf_{s\to \infty}F_r(y(s))\bigg)\\
&\leq&\frac{m}{\epsilon}.
\end{eqnarray*}
We now set
$$K:=\bigg\{y \in M \times V^*: 
\mathrm{dist}_{PS}(y,K_0)\leq \frac{m}{\epsilon}
\bigg\}.$$
Since $g_{PS}$ is geodesically complete, the set $K$ is compact. Moreover,
the estimate above shows that
$$\Bigg\{y(\sigma):y \in \bigcup_{r \in [0,1]} 
  \mathcal{M}(F_r,g_r),\,\,\sigma \in \mathbb{R}\Bigg\}
  \subset K$$
holds. This proves Step 6. 
\\ \\
\textbf{Step 7: }\emph{We prove the theorem}
\\ \\
It remains to show that the functions $F_r$ are Morse-Bott. We prove that in
Proposition~\ref{lagmorbo} below. This finishes the proof of the theorem.
\hfill $\square$
\\

If $x$ is a critical point of a Morse-Bott function, then we define the
index $\mathrm{ind}_F(x)$ of $F$ at $x$ as the number of negative eigenvalues
of the Hessian of $F$ at $x$. Note that the Morse-Bott condition implies
that the index is constant on each connected component of 
$\mathrm{crit}(F)$. The following proposition shows that if 
$f|_{h^{-1}(0)}$ is Morse-Bott, then the Lagrange multiplier functional is
also Morse-Bott and its index is independent of the behaviour of
$f$ outside of $h^{-1}(0)$. However note, that the Hessian itself depends
also on the derivatives of $f$ in the normal directions of $h^{-1}(0)$.

\begin{prop}\label{lagmorbo}
Let $M$ be a finite dimensional manifold and let $(V,\langle\,,\,\rangle)$ 
be a $k$-dimensional Euclidean vector space. Assume that  
$f \in C^\infty(M)$, 
$h \in C^\infty(M,V)$ such that $0$ is a regular value of $h$ and
$f|_{h^{-1}(0)}$ is Morse-Bott. Then $F \in C^\infty(M\times V^*)$
defined by
$F(x,v^*)=f(x)+v^*(h(x))$ for $(x,v^*) \in M \times V^*$
is also Morse-Bott. Moreover, if $\lambda \colon \mathrm{crit}(F)
\to \mathrm{crit}(f|_{h^{-1}(0)})$ is the natural bijection given by
$(x,v^*) \mapsto x$ for $(x,v^*) \in \mathrm{crit}(F)$,
then for the indices the following relation holds
$$\mathrm{ind}_F(\lambda^{-1}(x))=\mathrm{ind}_{f|_{h^{-1}(0)}}(x)+k,
\quad x \in \mathrm{crit}(f|_{h^{-1}(0)}).$$
\end{prop}
\textbf{Proof:} Let $x \in \mathrm{crit}(f|_{h^{-1}(0)})$. 
We first choose convenient coordinates around
$x$ in $M$. Set $n=\mathrm{dim}(M)$ and chose
$\delta_1,\delta_2>0$ so small
such that there exists a diffeomorphism $\phi$ from
$B_{\delta_1}^{n-k} \times B_{\delta_2}^k=
\{v \in \mathbb{R}^{n-k}: ||v||<\delta_1\}\times
\{v \in \mathbb{R}^k: ||v||<\delta_2\}$ to an open neighbourhood 
$U$ of $x$ in $M$ such that
$$\phi(0,0)=x, \quad h(\phi(q,w))=w, \quad q \in B_{\delta_1}^{n-k},
\,\,w \in B_{\delta_2}^k.$$
Choose furthermore an orthonormal basis in $V^*$ to define 
an isomorphism $\Phi \colon V^* \to \mathbb{R}^k$.
Let $\hat{f} \in C^\infty(B_{\delta_1}^{n-k}\times B_{\delta_2}^k)$,
be given by the pullback of $f$, i.e.
$$\hat{f}=\phi^* f|_{U},$$
and 
$\hat{F} \in C^\infty(B_{\delta_1}^{n-k}\times B_{\delta_2}^k\times
\mathbb{R}^k)$
be given by the pullback of $F$, i.e.
$$\hat{F}=(\phi \times \Phi)^*F|_{U \times V^*}.$$
Then $\hat{F}$ reads
$$\hat{F}(q,w,v)=\hat{f}(q,v)+\langle v,w \rangle,
\quad q \in B_{\delta_1}^{n-k},\,\,v \in B_{\delta_2}^k,\,\,
w \in \mathbb{R}^k.$$
We choose the standard flat metric on
$B_{\delta_1}^{n-k}\times B_{\delta_2}^k \times \mathbb{R}^k$
and introduce the $k\times k$-matrix $A$, the $k \times (n-k)$-matrix $B$
and the $(n-k) \times (n-k)$-matrix $H$ by
$$A_{ij}=\frac{\partial^2 \hat{f}(0,0)}{\partial w_i \partial w_j},
\quad
B_{ij}=\frac{\partial^2 \hat{f}(0,0)}{\partial q_i \partial w_j},
\quad
H_{ij}=\frac{\partial^2 \hat{f}(0,0)}{\partial q_i \partial q_j}.$$
Denote by $\pi_2 \colon M \times V^* \to V^*$ the
projection to the second factor.
The Hessian of 
$\hat{F}$ at $(0,0,\Phi \circ \pi_2 \circ \lambda^{-1}(x))$ 
with respect to the standard flat metric
is given by
$$H_{\hat{F}}(0,0,\Phi \circ \pi_2 \circ \lambda^{-1}(x))=
\left(\begin{array}{ccc}
H & B & 0\\
B^T & A & \mathrm{id}\\
0 &\mathrm{id} & 0
\end{array}
\right)$$
 We claim that
\begin{equation}\label{hessclaim}
\mathrm{dim}\big
(\mathrm{ker} H_{\hat{F}}(0,0,\Phi\circ \pi_2 \circ \lambda^{-1}(x))\big)
=\mathrm{dim}\big(\mathrm{ker} H\big).
\end{equation}
To see that assume that the vector $(\hat{q},\hat{w},\hat{v})
\in \mathbb{R}^{n-k}\times \mathbb{R}^k \times \mathbb{R}^k$
lies in the kernel of $H_{\hat{F}}(0,0,\Phi \circ \pi_2 \circ \lambda^{-1}(x)$.
It follows that
$$
\left\{\begin{array}{c}
H \hat{q}+B\hat{w}=0\\
B^T \hat{q}+A \hat{w}+\hat{v}=0\\
\hat{w}=0
\end{array}\right.$$
which implies that
$$(\hat{q},\hat{w},\hat{v})=(\hat{q},0,-B^T \hat{q}),\quad \hat{q} \in
\mathrm{ker}(H).$$
Hence (\ref{hessclaim}) follows.
\\
To prove that $F$ is Morse-Bott we denote for $y \in \mathrm{crit}(F)$ by
$\mathrm{dim}_y (\mathrm{crit}(F))$ the local dimension at $y$ of the
(unconnected) manifold $\mathrm{crit}(F)$ and compute
using (\ref{hessclaim}) and the Morse-Bott assumption on $f|_{h^{-1}(0)}$ 
\begin{eqnarray*}
\mathrm{dim}\big(\mathrm{ker} H_{\hat{F}}(\lambda^{-1}(0)\big)
&=&\mathrm{dim}\big(\mathrm{ker}H_{f|_{h^{-1}(0)}}(x)\big)\\
&=&\mathrm{dim}_x(\mathrm{crit}(f|_{h^{-1}(0)}))\\
&=&\mathrm{dim}_{\lambda^{-1}(x)}(\mathrm{crit}(F)).
\end{eqnarray*}
This proves that $F$ is Morse-Bott.

It remains to compute the index of the Hessian of $F$. To do that we
consider the smooth family of functions 
$\hat{f}_r \in C^\infty(B_{\delta_1}^{n-k}\times B_{\delta_2}^k)$
for $r \in [0,1]$ defined by
$$\hat{f}_r(q,w)=(1-r)\hat{f}(q,w)+r\hat{f}(q,0),\quad
q \in B_{\delta_1}^{n-k},\,\,w \in B_{\delta_2}^k.$$ 
Then
$$\hat{f}_0=\hat{f}, \quad \hat{f}_r|_{B_{\delta_1}^{n-k}\times\{0\}}
=\hat{f}|_{B_{\delta_1}^{n-k}\times \{0\}},\,\,
r \in [0,1].$$
We define the smooth family of functions
$\hat{F}_r \in C^\infty(B_{\delta_1}^{n-k}\times B_{\delta_2}^k
\times \mathbb{R}^k)$ for $r \in [0,1]$ by
$$\hat{F}_r(q,w,v)=\hat{f}_r(q,w)+\langle v,w\rangle.$$
Define further the smooth family of vectors $v_r \in \mathbb{R}^k$
for $r \in [0,1]$ by
$$(v_r)_i=-\frac{\partial \hat{f}_r(0,0)}{\partial w_i},\quad
i \in \{1,\ldots,k\}.$$
The functions $\hat{F}_r$ have critical points at $(0,0,v_r)$
and it follows from (\ref{hessclaim}) that 
$\mathrm{dim}\big(\mathrm{ker} H_{\hat{F}_r}(0,0,v_r)\big)=
\mathrm{dim}\big(\mathrm{ker}(H)\big)$ does not depend on
$r \in [0,1]$. Since the eigenvalues of a continuous family of 
matrices are continuous, see \cite[Theorem II.5.1]{kato} we conclude
\begin{equation}\label{hess1}
\mathrm{ind}_{\hat{F}_0}(0,0,v_0)
=\mathrm{ind}_{\hat{F}_1}(0,0,v_1).
\end{equation}
The Hessian of $\hat{F}_1$ at $(0,0,v_1)=(0,0,0)$
is given by
$$H_{\hat{F}_1}(0,0,0)=
\left(\begin{array}{ccc}
H & 0 & 0\\
0 & 0 & \mathrm{id}\\
0 &\mathrm{id} & 0
\end{array}
\right)$$
from which we deduce
\begin{equation}\label{hess2}
\mathrm{ind}_{\hat{F}_1}(0,0,v_1)=
\mathrm{ind}_{\hat{f}|_{B_{\delta_1}^{n-k}\times\{0\}}}(0,0)+k.
\end{equation}
Combining (\ref{hess1}) and (\ref{hess2}) we compute
\begin{eqnarray*}
\mathrm{ind}_F(\lambda^{-1}(x))&=&
\mathrm{ind}_{\hat{F}_0}(0,0,v_0)\\
&=&\mathrm{ind}_{\hat{f}|_{B_{\delta_1}^{n-k}\times\{0\}}}(0,0)+k\\
&=&\mathrm{ind}_{f|_{h^{-1}(0)}}(x)+k.
\end{eqnarray*}
This completes the proof of the proposition. \hfill $\square$

\end{document}